\documentclass[12pt,psfig]{article}
\usepackage{graphicx,epsfig}
\newcommand{\ti}{\tilde}

\newcommand{\ba}{\begin{array}}
\newcommand{\ea}{\end{array}}
\newcommand{\ban}{\begin{eqnarray*}}
\newcommand{\ean}{\end{eqnarray*}}

\newcommand{\ol}{\overline}
\newcommand{\lf}{\left}
\newcommand{\rf}{\right}

\newcommand{\la}{\langle}
\newcommand{\ra}{\rangle}

\begin{document}

\baselineskip=17pt

\begin{center}

\vspace{-0.6in}
{\large \bf  Boundedness and closedness of linear relations$^{\dag}$}

 \vspace{0.3in}

Yuming Shi$^{a,{\ddag}}$, Guixin Xu$^a$, Guojing Ren$^b$ \\

$^a$Department of Mathematics, Shandong University\\

Jinan, Shandong 250100, P. R. China\\

$^b$School of Mathematics and Quantitative Economics\\

Shandong University of Finance and Economics\\

Jinan, Shandong 250014, P. R. China\\

\footnote{$^\dag$This research was supported by the NNSF of China (Grants 11571202
and 11301304).\\
\indent \ \ $^\ddag$The corresponding author.\\
\indent \indent Email addresses: ymshi@sdu.edu.cn(Y. Shi), guixinxu$_-$ds@163.com(G. Xu), gjren@sdufe.edu.cn(G. Ren).}

\end{center}

{\bf Abstract.} This paper studies boundedness and closedness of linear relations,
which include both single-valued and multi-valued linear operators.
A new (single-valued) linear operator induced by a linear relation
is introduced, and its relationships with other two important
induced linear operators  are established.
Several characterizations for closedness, closability,
bundedness, relative boundedness, and boundedness from below (above)
of linear relations are given in terms of their induced linear operators.
In particular, the closed graph theorem for linear relations in Banach spaces
is completed, and stability of closedness of
linear relations under bounded and relatively bounded perturbations is studied.
The results obtained in the present paper generalize
the corresponding results for (single-valued) linear operators to multi-valued linear operators,
and some improve or relax certain assumptions of the related existing results.

\medskip

\noindent{\bf 2010 AMS Classification}: 47A06, 47B15, 47A55.

\bigskip

\noindent{\bf 1. Introduction}\medskip

In the study of classical linear operators, it is always required that
 the operators are single-valued. They have been extensively studied
 and a great deal of elegant results have been obtained (cf., [10, 14, 25, 29, 36]).

Recently, it was found that minimal and maximal operators generated by symmetric
linear expressions in the discrete and time scales cases are multi-valued or non-densely defined
in general even though the corresponding definiteness condition is satisfied (cf.,
[26, 34]), and similar are those generated by symmetric linear differential expressions that do not
satisfy the definiteness condition [19].
So we should apply the theory of multi-valued linear operators to
study the above problems instead of the classical operator theory.
This is the motivation for us to study some topics about multi-valued
linear operators, which are a necessary foundation of research on
those related problems about difference operators
  as well as operators defined on time scales.

A multi-valued linear operator from a linear space $X$ to a linear
space $Y$ is also called a linear relation.
Since its graph is a linear subspace (briefly, subspace) of the product space $X\times Y$,
it is also called a subspace [7]. It is more convenient to study
many problems about multi-valued linear operators by their graphs.
So we shall use the term ``subspace" in the present paper.
Throughout the present paper, a linear operator always means
a single-valued linear operator for convenience.

To the best of our knowledge, the theory of subspaces of product spaces
 was initiated by von Neumann [23, 24].
The operational calculus of subspaces was developed by Arens [3].
Their works were followed by many scholars,
 and some basic concepts, fundamental properties, extension, resolvent,
 spectrum, and perturbation for subspaces were studied (cf., [1-2, 4-9, 12-13, 17-18, 30-33, 37]).
The theory of subspaces has been successfully applied in the analysis of linear and
nonlinear problems, control theory, and linear difference
equations (cf., [15-16, 19-22, 26-28, 34-35]).

The multi-valued part of a subspace
will result in the main difficulty in the study of subspaces.
In order to deal with it, some scholars introduced different
``operator parts" of subspaces, which provide a bridge between
subspaces and linear operators, so
that one can apply the theory of linear operators to study
some properties of subspaces. This term was coined by Coddington [8].
In 1961, Arens decomposed a closed subspace $T$ in $X^2$ as an orthogonal sum of
a singled-valued  operator part $T_s$ and a purely multi-valued part $T_{\infty}$,
where $X$ is a Hilbert space [3].
This decomposition has been well applied in our study of subspaces (cf., [30-33]).
Note that it is required in this decomposition
that the space $X$ is a Hilbert space and the subspace $T$ is closed.
In 1990 and 1991, Lee and Nashed introduced algebraic operator part (also called algebraic selection)
for a subspace of the product space of two linear spaces, and topological and proximinal
operator parts (also called topological and proximinal selections) for a subspace
of the product space of two normed linear spaces [17, 18].
In 1998, Cross defined a linear operator, denoted by ${\tilde T_s}$, through multiplying a
related natural quotient map for any subspace $T$ of the product space
of two linear spaces [9].
Note that $T_s$, and the algebraic, topological, and proximinal
selections of a subspace $T$ are subspaces of $T$, but ${\tilde T_s}$
is not a part of $T$. However, the operator ${\tilde T_s}$ is very convenient
in the study of some problems. We shall introduce a new linear operator induced by
a subspace in the present paper.

There are still many important fundamental problems about subspaces
that have neither been studied nor completed.
It is well known that the  closed
and bounded operators are very important two classes in the theory of linear operators.
The closedness and boundedness of subspaces have been, but not thoroughly studied.
In the present paper, we shall focus on these two classes of subspaces.

Now, we shall recall some existing research works related to the boundedness and closedness
of subspaces. The concept of closedness for a subspace is defined by its closedness
in the corresponding product space. This is the same as that for a linear operator by its graph,
which is a subspace of the product space. In [9], Cross extended the concept
of closability for a linear operator to a subspace and gave some characterizations
for closedness and closability of subspaces.
Since a subspace is multi-valued in general,
it is not easy to introduce its norm  at one point and its norm.
In 1972, Robinson was motivated by some problems of convex analysis
 and mathematical programming,
and defined a norm of a convex process on Banach spaces [28]. In 1991, Lee and Nashed
adopted this definition for subspaces, and studied the relations between the norm of a
subspace and the norms of its algebraic operator parts [17]. They gave
some sufficient conditions under which the infimum of the norms of its algebraic operator parts
is attained. They pointed out that for a given subspace $T$ with finite norm, there
may not exist any algebraic operator part $R$ such that $\|R\|=\|T\|$ in general.
In 1998, Cross defined a norm of a subspace $T$ at one point and its norm
by ${\tilde T_s}$ [9], and shown that this norm
is equal to the norm given by Lee and Nashed in [17]. Concepts
of boundedness and relative boundedness of subspace
can be defined by their norm and their norms at points, respectively.

In the present paper,
we shall give some new equivalent characterizations for
closedness and closability of subspaces, equivalent characterizations and sufficient conditions
for bundedness and relative boundedness of subspaces,
and equivalent characterizations for boundedness from below (above), non-negativeness
(non-positiveness), and positiveness (negativeness) of Hermitian subspaces.
In particular, we shall study the stability of closedness of
subspaces under bounded and relatively bounded perturbations.

The rest of this paper is organized as follows.
In Section 2, some notations, basic concepts,
and fundamental results about closed, closable, and bounded
subspaces are introduced. Three operators induced by a subspace
are introduced and their relationships are studied.
In Section 3, some equivalent characterizations
for boundedness and relative boundedness of subspaces,
and boundedness from below (above),
non-negativeness (non-positiveness), and positiveness (negativeness) of
Hermitian subspaces are given in terms of their induced operators, respectively.
Some properties and characterizations
for closedness and closability of subspaces, and the closed
graph theorem for subspaces are studied in Sections 4 and 5,  respectively.
Finally, two sufficient conditions
for relative boundedness, and the stability of the closedness
of subspaces under bounded and relatively bounded perturbations
are discussed in Section 6.

\bigskip

\noindent{\bf 2. Preliminaries}\medskip

In this section, we shall introduce some notations and basic concepts,
including closed, closable, and bounded subspaces, and
give some fundamental results about subspaces. In particular,
we shall introduce three operators induced by a subspace, which one is new,
and discuss their relationships. These operators will play an important role in the study
of the theory of subspaces.

This section is divided into three subsections.

\bigskip

\noindent{\bf 2.1. Some notations and basic concepts about subspaces}\medskip


Let $X, Y$, and $Z$ be linear spaces over a number field ${\mathbf K}$.
If $X$ is a normed space with norm $\|\cdot\|_X$ or an inner product space with inner
product $\langle \cdot, \cdot\rangle_X$,
the subscript $X$ will be omitted without confusion.
If $X$ is an inner product space and $E\subset X$, by $E^\perp$
denote the orthogonal complement of $E$.

In the case that $X$ and $Y$ are topological linear spaces,
the topology of the product space $X\times Y$ is naturally induced by $X$ and $Y$.
If $X$ and $Y$ are normed,
then the norm of $X\times Y$ is defined
by
\vspace{-0.2cm}$$\|(x,y)\|=\left(\|x\|^2+\|y\|^2\right)^{1/2},\;(x,y)\in X\times Y.\vspace{-0.2cm}$$
Similarly, if $X$ and $Y$ are inner product spaces,  then the inner product
of $X\times Y$ is defined by
\vspace{-0.2cm}$$\langle(x_1,y_1),(x_2,y_2)\rangle=\langle x_1, x_2\rangle
+\langle y_1,y_2\rangle,\;\;(x_1,y_1),\;(x_2,y_2)\in X\times Y.$$

By $LR(X,Y)$ denote
the set of all the linear subspaces (briefly, subspaces)
of the product space $X\times Y$. In the case that $Y=X$, briefly by $LR(X)$ denote
$LR(X,Y)$.

Let $T\in LR(X,Y)$. By $D(T)$ and $R(T)$ denote the domain and range of $T$,
respectively. Further, denote
\vspace{-0.2cm}$$T(x):=\{ y\in Y:\; (x, y)\in T\}, \; N(T):=\{x\in X:\;(x,0)\in T\},\vspace{-0.4cm}$$
\vspace{-0.4cm}$$T^{-1}:=\{ (y,x):\; (x, y)\in T\}.\vspace{-0.1cm}$$
It is evident that $T(0) = \{0\}$ if and only if
$T$ uniquely determines a linear operator
from $D(T)$ into $Y$ whose graph is $T$. For convenience,
a linear operator from $X$ to $Y$ will always be
identified with a subspace of $X\times Y$ via its graph.

Let $T, S, W\in LR(X,Y)$ and $\alpha\in {\mathbf K}$. Define
\vspace{-0.2cm}$$\begin{array}{ccc}\alpha\, T :
= \{(x, \alpha\, y) : (x, y) \in T\},\\[0.5ex]
T + S := \{(x, y + z) : (x, y) \in T,\; (x, z)\in S\}.
\end{array}\vspace{-0.2cm}$$
It can be easily verified that the above sum satisfies the laws of commutation and association:
\vspace{-0.2cm}$$T+S=S+T,\;\; (T+S)+W=T+(S+W).\vspace{-0.1cm}$$
On the other hand, if $T \cap S$ = \{(0, 0)\}, then denote
\vspace{-0.2cm}$$T\dot{+} S:=\{(x+u,y+v): \;(x,y)\in T,\; (u,v)\in S\}.\vspace{-0.2cm}$$
Further, in the case that $X$ and $Y$ are inner product spaces,
if $T$ and $S$ are orthogonal; that is,
$\langle(x, y),(u,v)\rangle = 0 $ for all
$(x, y) \in T$ and $(u, v) \in S$, then denote
\vspace{-0.2cm}$$T\oplus S:=T\dot{+} S.\vspace{-0.2cm}$$
Let $T\in LR(X,Y)$ and $S\in LR(Y,Z)$. The product of $T$ and $S$
is defined by (see [3])
\vspace{-0.2cm}$$ST=\{(x,z)\in X\times Z:\;{\rm there\; exists}\; y \in Y\;
{\rm such\; that}\;(x,y)\in T\;{\rm and}\; (y,z)\in S\}.                             \vspace{-0.2cm}$$
Note that if $S$ and $T$ are operators, then $ST$ is also an operator.

If $X$ and $Y$ are topological linear spaces,
then, by $CR(X,Y)$ denote the set of all the closed subspaces of $X\times Y$.
In the case that $X=Y$, by $CR(X)$ denote $CR(X,Y)$ briefly.
Let $T\in LR(X,Y)$. By ${\ol T}$ denote the closure of $T$.
Obviously, ${\ol T}\in CR(X,Y)$. Subspace $T$ is said to be closed
if ${\ol T}=T$, and closable if ${\ol T}(x)=T(x)$ for each $x\in D(T)$. It is evident that
$T\in CR(X,Y)$ if and only if $T^{-1}\in CR(Y,X)$.

\medskip

\noindent{\bf Lemma 2.1 {\rm [9]}.} Let $X$ and $Y$ be linear spaces
and $T\in LR(X,Y)$. Then $T(x)=\{y\}+T(0)$ for
every $x\in D(T)$ and $y\in T(x)$.

\bigskip

\noindent{\bf 2.2. Three operators induced by a subspace
and their relationships}\medskip

In this subsection, we shall first recall two linear operators
induced by a subspace, given by Cross and Arens [3, 9], separately,
and then introduce another new linear operator induced by the subspace.
Finally, we shall establish their relationships.

We shall first introduce a natural quotient map and then define an operator induced
by a subspace, which was gave in [9, Section II.1]. Let $X$ be a linear space,
and $E$ be a subspace of $X$. Define the following quotient space [14]:
\vspace{-0.2cm}$$X/E:=\{[x]:\;x\in X\},\;[x]:=\{x\}+E.\vspace{-0.2cm}$$
If $X$ is normed and $E$ is closed, then $X/E$  is a normed linear space with norm
\vspace{-0.2cm}$$\|[x]\|:=d(x,E)=\inf\{\|x-e\|:\;e\in E\},\;[x]\in X/E.\vspace{-0.2cm}$$
Further, if $X$ is complete, so is $X/E$.
Let $X$ be a Hilbert space and $E$ be a closed subspace of $X$.
We define
an inner product on the quotient space $X/E$ by
\vspace{-0.2cm}$$\la [x],[y]\ra=\la x^\perp,y^\perp\ra,\;\;[x],[y]\in X/E,              \eqno (2.1)\vspace{-0.2cm}$$
where $x=x_0+x^\perp$, $y=y_0+y^\perp$ with $x_0,y_0\in E$ and $x^\perp,y^\perp\in E^\perp.$
It can be easily verified that the above inner product is well-defined and
$X/E$ with this inner product is a Hilbert space. The norm induced
by this inner product is the same as that of $X/E$ induced by the norm
of $X$.

Now, define the following natural quotient map:
\vspace{-0.2cm}$$Q^X_E:\;X\to X/E,\;\;x\mapsto [x].\vspace{-0.2cm}$$
Let $T\in LR(X,Y)$ and $Y$ be a topological linear space. By $Q_T$ denote $Q^Y_{\overline{T(0)}}$
for briefness without confusion.
Define
\vspace{-0.2cm}$${\tilde T_s}:=G(Q_T)\, T.                                            \eqno (2.2)\vspace{-0.2cm}$$
Then ${\ti T_s}\in LR(X, Y/\,{\overline{T(0)}})$ is a linear operator
with domain $D(T)$ [9, Proposition II.1.2].

\medskip

\noindent{\bf Lemma 2.2.} Let $X$ and $Y$ be topological
linear spaces and $T\in LR(X,Y)$. If $T$ is closed, then $T(0)$ is closed.
Further, if $Y$ is a normed space, then $T$ is closed if and only if
${\ti T_s}$ and $T(0)$ are closed.\medskip

\noindent{\bf Proof.} The first assertion can be directly obtained by the
definition of closedness of a subspace.

Further, suppose that $Y$ is a normed space. With a similar argument
to that used in the proof of [9, Proposition II.5.3],
one can show that the second assertion holds. The proof is complete.
\medskip

\noindent{\bf Lemma 2.3.}
Let $X$ and $Y$ be topological linear spaces and $T\in LR(X,Y)$. Then
${\ol {T(0)}}\subset {\ol T}(0).$
\medskip

\noindent{\bf Proof.} Since ${\ol T}$ is closed, ${\ol T}(0)$ is closed by Lemma 2.1.
Hence, the assertion holds by the fact that $T(0)\subset {\ol T}(0)$,
and the proof is complete.\medskip

\noindent{\bf Remark 2.1.} Cross gave the same results
as the second result of Lemma 2.2 and that of Lemma 2.3
in the case that $X$ and $Y$ are both normed spaces in [9, Proposition II.5.3
and Exercise II.5.19].\medskip

In the case that $X$ is a Hilbert space and $T\in CR(X)$, Arens introduced another operator
by the following decomposition [3]:
\vspace{-0.2cm}$$T=T_s\oplus T_\infty,                                                                \vspace{-0.2cm}$$
where
\vspace{-0.2cm}$$T_\infty :=\{(0,g)\in X^2: (0,g)\in T\},\;\;
T_s :=T\ominus T_{\infty}.\vspace{-0.1cm}$$
Then $T_s\in CR(X)$ is a linear operator
with domain $D(T)$, and $T_\infty\in CR(X)$.
So $T_s$ and $T_\infty$ are often called the operator
and pure multi-valued parts of $T$, respectively.
In addition,
they satisfy the following properties [3]:
\vspace{-0.2cm}$$D(T_s)=D(T),\;R(T_s)\subset T(0)^\perp,\;
\;T_\infty =\{0\}\times T(0).                                                       \eqno (2.3)\vspace{-0.2cm}$$

Note that $T_s$ is a subspace of $T$,
but ${\ti T_s}$ is not. Before giving out the relationship
between ${\ti T_s}$ and $T_s$, we introduce another
new operator for $T$ as follows.\medskip

\noindent{\bf Lemma 2.4.} Let $X$ be a Hilbert space
and $E$ be a closed subspace of $X$.
Then the map
\vspace{-0.2cm}$$Q_E^X|_{E^\perp}:\;E^\perp\to X/E,\;\;x\mapsto [x],\vspace{-0.2cm}$$
is isomorphic, and preserves norm. \medskip

\noindent{\bf Proof.} It can be easily verified that $Q_E^X|_{E^\perp}$ is linear and bijective.
Further, it satisfies that
\vspace{-0.2cm}$$\|Q_E^X|_{E^\perp}(x)\|=\|[x]\|=\|x\|,\;\,x\in E^\perp.               \eqno (2.4)\vspace{-0.2cm}$$
Hence, the assertion holds. The proof is complete.
\medskip

Let $X$ be a linear space, $Y$ be a Hilbert space, and $T\in LR(X,Y)$. By $P_T$
denote the following orthogonal projection:
\vspace{-0.2cm}$$P_T:\, Y\to {\overline {T(0)}}^\perp.\vspace{-0.2cm}$$
Then
\vspace{-0.2cm}$$Q_T=Q_T|_{{\overline {T(0)}}^\perp}\,P_T,\;\;
P_T=\left(Q_T|_{{\overline {T(0)}}^\perp}\right)^{-1}Q_T.                        \eqno (2.5)\vspace{-0.2cm}$$
Define
\vspace{-0.2cm}$${\hat T_s}:=G(P_T)\,T.                                              \eqno (2.6)   \vspace{-0.2cm}$$
Then ${\hat T_s}\in LR(X,Y)$ with $D({\hat T_s})=D(T)$ and $R({\hat T_s})
\subset {\overline {T(0)}}^\perp$. Note that ${\hat T_s}=T$
in the case that $T$ is single-valued.

\medskip

\noindent{\bf Theorem 2.1.} Let $X$ be a linear space, $Y$ be a Hilbert space,
and $T\in LR(X,Y)$.
Then
\vspace{-0.2cm}$${\ti T_s}=G(Q_T|_{{\ol {T(0)}}^\perp})\,{\hat T_s}.                   \eqno (2.7)\vspace{-0.2cm}$$
And consequently, ${\hat T_s}$ is a linear operator with domain $D(T)$.
Further, if $X$ is a topological linear space,
then ${\hat T_s}\subset {\ol T}$, which implies that
${\hat T_s}\subset T$ in the case that $T$ is closed.\medskip

\noindent{\bf Proof.}  Assertion (2.7) can be easily derived from (2.2), (2.5),
(2.6), and the fact that $G(W)=G(V)G(U)$ if $W=VU$ for linear operators $U$ and $V$.
Hence, ${\hat T_s}$ is a linear operator with domain $D(T)$ by the fact that
${\ti T_s}$ is single-valued.

Further, suppose that $X$ is a topological linear space.
For any fixed $(x, z)\in {\hat T_s}$, there exist
$y\in Y$ and $y_1\in {\ol {T(0)}}$ such that $(x,y)\in T$ with $y=z+y_1$ by (2.6).
In addition, it follows from Lemma 2.3 that ${\ol {T(0)}}\subset {\ol T}(0)$,
which implies that $y_1\in {\ol T}(0)$. Hence,
$(x,z)=(x,y)-(0,y_1)\in {\ol T}$. Therefore, ${\hat T_s}\subset {\ol T}$.
And consequently, ${\hat T_s}\subset T$ if $T$ is closed.
The proof is complete.\medskip

In the case that $X$ is a topological linear space, $Y$ is a Hilbert space,
 and $T\in CR(X,Y)$, ${\hat T_s}$ can be called an operator part of $T$.
The following result gives the relationship among three operators
${\ti T_s}$, $T_s$, and ${\hat T_s}$.\medskip

\noindent{\bf Theorem 2.2.} Let $X$ be a Hilbert space and $T\in CR(X)$.
Then
\vspace{-0.2cm}$$T_s={\hat T_s}=\left(G(Q_T|_{T(0)^\perp})\right)^{-1}{\ti T_s}.                \eqno (2.8)$$

\noindent{\bf Proof.} Since $T$ is closed, $T(0)$ is closed by Lemma 2.2.
Then the first equality in (2.8) directly follows from the definitions of $T_s$
and ${\hat T_s}$, and the second equality in  (2.8) is derived from (2.7).
This completes the proof.\medskip

Since ${\ti T_s}(x)$, ${\hat T_s}(x)$,
and $T_s(x)$ contain singleton elements, respectively,
for every $x\in D(T)$, by ${\ti T_s}(x)$, ${\hat T_s}(x)$,
and $T_s(x)$ also denote their elements
as usual operators for convenience when it is needed
in the rest of this paper.

\bigskip

\noindent{\bf 2.3. Properties of norms of subspaces}\medskip

In this subsection, we shall introduce concepts of the norm of a subspace
 at one point in its domain, the norm of a subspace, and a bounded subspace,
and discuss their fundamental properties.

Let $X$ and $Y$ be normed spaces and $T\in LR(X,Y)$. The norm of $T$ at $x\in D(T)$
and the norm of $T$ are defined by, respectively (see [9, Section II.1]),
\vspace{-0.2cm}$$\|T(x)\|:=\|{\ti T_s}(x)\|,\;\;
\|T\|:=\|{\ti T_s}\|=\sup\{\|{\ti T_s}(x)\|:\,x\in D(T)\;{\rm with}\;\|x\|\le 1\}.             \eqno (2.9)$$

\noindent{\bf Definition 2.1.} Let $X$ and $Y$ be normed spaces and $T\in LR(X,Y)$.
If $\|T\|$ is finite, then $T$ is said to be bounded.
\medskip

By $BR(X,Y)$ denote the set of all the bounded subspaces of $X\times Y$. Note that
if $T$ is an operator, then $T$ is bounded in the sense of operator
if and only if its graph $G(T)$ is bounded in the sense of subspace with the same bound.
\medskip

\noindent{\bf Remark 2.2.} In [9], Cross gave another definition of
 a bounded subspace $T$ (see [9, Definition II.1.3]), in which it is required that
$D(T)=X$. Now, we remove this requirement in Definition 2.1. This agrees with
that of a bounded operator [36].\medskip

Next, we recall some fundamental results about norms of subspaces.\medskip

\noindent{\bf Lemma 2.5 {\rm [9]}.}
Let $X$ and $Y$ be normed spaces over a field ${\mathbf K}$,
and $S,T\in LR(X,Y)$. Then
\begin{itemize}\vspace{-0.2cm}
\item[{\rm (i)}]$\|T(x)\|=d(y,{\ol {T(0)}})=d(y,T(0))=d(0,T(x))=d(T(0),T(x))$
for every $x\in D(T)$ and every $y\in T(x)$;
\vspace{-0.2cm}
\item[{\rm (ii)}]$\|(\alpha\, T)(x)\|=|\alpha|\,\|T(x)\|$ for every $x\in D(T)$
and every $\alpha\in {\mathbf K}$;
\vspace{-0.2cm}
\item[{\rm (iii)}]$\|(S+T)(x)\|\le \|S(x)\|+\|T(x)\|$ for every $x\in D(T)$.

\end{itemize}

\medskip

\noindent{\bf Lemma 2.6 {\rm [9, 18]}.}
Let $X$ and $Y$ be normed spaces over a field ${\mathbf K}$, and
$S,T\in LR(X,Y)$. Then
\begin{itemize}\vspace{-0.2cm}
\item[{\rm (i)}]$\|T\|=\sup\{\|T(x)\|:\;x\in D(T)\;{\rm with}\;\|x\|\le 1\}$;
\vspace{-0.2cm}
\item[{\rm (ii)}]$\|T(x)\|\le \|T\|\|x\|$ for all $x\in D(T)$ in the case that $T$ is bounded;

\vspace{-0.2cm}
\item[{\rm (iii)}]$\|\alpha\, T\|=|\alpha|\,\|T\|$ for every $\alpha\in {\mathbf K}$;
\vspace{-0.2cm}
\item[{\rm (iv)}]$\|S+T\|\le \|S\|+\|T\|$.

\end{itemize}

Note that the norm $\|T\|$ is not a real norm since the following inequalities
may not hold in general (see [9, Exercise II.1.12]):
\vspace{-0.2cm}$$\|(S-T)(x)\|\ge \|S(x)\|-\|T(x)\|\;{\rm for}\;x\in D(S)\cap D(T),
\;\|S-T\|\ge \|S\|-\|T\|.\vspace{-0.2cm}$$
We shall show that the above inequalities hold under some conditions.\medskip

\noindent{\bf Proposition 2.1.} Let $X$ and $Y$ be linear spaces, and $S, T\in LR(X,Y)$.
Then $S=(S-T)+T$  if and only if $D(S)\subset D(T)$ and $T(0)\subset S(0)$.
\medskip

\noindent{\bf Proof.} ``$\Rightarrow$"  Suppose that $S=(S-T)+T$.
Since $D((S-T)+T)=D(S)\cap D(T)$, we have that $D(S)\subset D(T)$.
In addition, $S(0)=(S-T)(0)+T(0)$ and $0\in (S-T)(0)$.
It follows that $T(0)\subset S(0)$.

``$\Leftarrow$" Suppose that $D(S)\subset D(T)$ and $T(0)\subset S(0)$.
Then $D((S-T)+T)=D(S)\cap D(T)=D(S)$. Further, for any $x\in D(S)$, we have that
\vspace{-0.2cm}$$\left((S-T)+T\right)(x)=(S-T)(x)+T(x)=S(x)-T(x)+T(x).                 \eqno (2.10) \vspace{-0.2cm}$$
Fix any $y\in S(x)$ and any $z\in T(x)$. By Lemma 2.1 we get that
\vspace{-0.2cm}$$S(x)=\{y\}+S(0),\;\;T(x)=\{z\}+T(0).                                  \eqno (2.11)\vspace{-0.2cm}$$
So it follows from (2.10) and (2.11) that $-T(x)+T(x)=T(0)$ and
\vspace{-0.2cm}$$\left((S-T)+T\right)(x)=\{y\}+S(0)+T(0)=\{y\}+S(0)=S(x).\vspace{-0.2cm}$$
Hence, $S=(S-T)+T$. This completes the proof.

\medskip

\noindent{\bf Remark 2.3.} The sufficiency of Proposition 2.1 was given in [2, (i) of Lemma 2.5]
in the case that $X$ and $Y$ are Banach spaces. In fact, it is only required
that $X$ and $Y$ are linear spaces in the proof of [2, (i) of Lemma 2.5].
Here, we give its detailed proof for completeness.

\medskip

\noindent{\bf Theorem 2.3.} Let $X$ and $Y$ be normed spaces,
and $S, T\in LR(X,Y)$
satisfy that $D(S)\subset D(T)$ and ${\ol {T(0)}}\subset {\ol {S(0)}}$. Then
\vspace{-0.2cm}$$\|(S-T)(x)\|\ge \|S(x)\|-\|T(x)\|,\;\;x\in D(S).                   \eqno (2.12)\vspace{-0.2cm}$$
Further, if $S$ and $T$ are bounded, then $S-T$ is bounded, and
\vspace{-0.2cm}$$\|S-T\|\ge \|S\|-\|T\|.                                               \eqno (2.13)\vspace{-0.2cm}$$
In addition, if either $T$ is bounded and $S$ is unbounded or
$S$ is bounded, $T|_{D(S)}$ is unbounded, and ${\ol {T(0)}}={\ol {S(0)}}$,
then $S-T$ is unbounded.\medskip

\noindent{\bf Proof.} It can be easily verified that ${\ol {S(0)}}={\ol {(S-T)(0)}}$
by the assumption that ${\ol {T(0)}}\subset {\ol {S(0)}}$. So by (i) of Lemma 2.5
we have that for any $x\in D(S)$, and any given $y_1\in S(x)$ and $y_2\in T(x)$,
\vspace{-0.2cm}$$\begin{array}{rrll}
&&\|(S-T)(x)\|=d(y_1-y_2,{\ol {(S-T)(0)}})=d(y_1-y_2,{\ol {S(0)}})\\[0.5ex]
&\geq & d(y_1,{\ol {S(0)}})-d(y_2,{\ol {S(0)}})
\geq  d(y_1,{\ol {S(0)}})-d(y_2,{\ol {T(0)}}) = \|S(x)\|-\|T(x)\|
\end{array}                                                                        \eqno (2.14)\vspace{-0.2cm}$$
which yields that (2.12) holds.

Further, suppose that $S$ and $T$ are bounded.
Then $S-T$ is bounded by (iii) of Lemma 2.5 and (i) of Lemma 2.6.
It follows from (2.14) and (ii) of Lemma 2.6 that
\vspace{-0.2cm}$$\|S(x)\|\le \|S-T\|+\|T\|\;{\rm for} \;x\in D(S)\;{\rm with}\; \|x\|\le 1. \vspace{-0.2cm}$$
Hence, $\|S\|\le \|S-T\|+\|T\|$ by (i) of Lemma 2.6. And consequently, (2.13) holds.

In addition, suppose that  $T$ is bounded and $S$ is unbounded.
It can be easily verified that $S-T$ is unbounded by (2.14) and (i) and (ii) of Lemma 2.6.

Finally, suppose that $S$ is bounded, $T|_{D(S)}$ is unbounded, and ${\ol {T(0)}}={\ol {S(0)}}$.
It follows from (2.14) that
\vspace{-0.2cm}$$\|(S-T)(x)\|\geq  d(y_2,{\ol {T(0)}})-d(y_1,{\ol {S(0)}})
= \|T(x)\|-\|S(x)\|,                                                  \vspace{-0.2cm}$$
which, together with (i) and (ii) of Lemma 2.6, implies that $S-T$ is unbounded.
The whole proof is complete.\medskip

\noindent{\bf Remark 2.4.} Let $X$ and $Y$ be normed spaces, and $S, T\in LR(X,Y)$
satisfy that $D(S)\subset D(T)$. If $S$ is bounded, $T$ is unbounded,
${\ol {T(0)}}\subset {\ol {S(0)}}$, and
${\ol {T(0)}}\neq {\ol {S(0)}}$, then $S-T$ may be bounded or unbounded.
For example, let $X=Y=l^2$, and
\vspace{-0.2cm}$$T(x)=\{\{nx(n)\}_{n=1}^\infty\},\;\;x=\{x(n)\}_{n=1}^\infty\in D(T),\vspace{-0.2cm}$$
where $D(T)=\{x=\{x(n)\}_{n=1}^\infty\in l^2:\,  \{nx(n)\}_{n=1}^\infty\in l^2\}$.
Then $T$ is unbounded and single-valued. Let $D(S_1)=D(S_2)=D(T)$,
and
\vspace{-0.2cm}$$S_1(x)=l^2,\;\;S_2(x)=\{x\}+S_2(0),\;\;x\in D(T),\vspace{-0.2cm}$$
where $S_2(0)={\rm span}\{e_1\}$ with $e_1(1)=1$ and $e_1(n)=0$ for $n\ge 2$.
It is evident that $S_1(0)=l^2$, $S_1$ is bounded with bound $\|S_1\|=0$,
and $S_2$ is bounded with bound $\|S_2\|=1$. In addition,
${\ol {T(0)}}=\{0\}\subset {\ol {S_i(0)}}$ and
${\ol {T(0)}}\neq {\ol {S_i(0)}}$ for $i=1,2$.
Further, we get that for any $x\in D(T)$,
\vspace{-0.2cm}$$S_1(x)-T(x)=l^2,\;\; S_2(x)-T(x)=\{\{(1-n)x(n)\}_{n=1}^\infty\}+S_2(0),\vspace{-0.2cm}$$
which implies that $S_1-T$ is bounded with bound $\|S_1-T\|=0$
and $S_2-T$ is unbounded.

\medskip

The following result gives the relationships among the norms of subspaces
$T$, ${\hat T_s}$, and $T_s$. It can be easily derived from (2.4), (2.7), and (2.8).\medskip

\noindent{\bf Theorem 2.4.} Let $X$ be a normed space, $Y$ be a Hilbert space,
and $T\in LR(X,Y)$. Then
\vspace{-0.2cm}$$\|{\hat T_s}(x)\|=\|T(x)\|\;{\rm for} \;x\in D(T);\;\;\|{\hat T_s}\|=\|T\|,        \eqno (2.15)  \vspace{-0.2cm}$$
and consequently, the boundedness of $T$ and ${\hat T_s}$ are equivalent.
Furthermore, if $X=Y$ is a Hilbert space and $T\in CR(X)$, then
\vspace{-0.2cm}$$\|T_s(x)\|=\|{\hat T_s}(x)\|=\|T(x)\|\;{\rm for} \; x\in D(T);\;\;
\|T_s\|=\|{\hat T_s}\|=\|T\|,                                                             \vspace{-0.2cm}$$
and consequently, the boundedness of $T$, ${\hat T_s}$, and $T_s$ are
equivalent in this special case.
\bigskip

\noindent{\bf 3. Boundedness and relative boundedness for subspaces}
\medskip

In this section, we shall give some sufficient and necessary conditions
for boundedness of subspaces, introduce concepts of
boundedness from below (above),  non-negativeness (non-positiveness), and positiveness
(negativeness) for Hermitian subspaces,
and relative boundedness for  subspaces,
and study their equivalent characterizations by their induced operators.

This section is divided into two subsections.
\bigskip

\noindent{\bf 3.1. Boundedness and boundedness from below (above) for subspaces}
\medskip

In this subsection, we shall first give some sufficient and necessary conditions
for boundedness of subspaces, then introduce concepts of
boundedness from below (above), non-negativeness (non-positiveness), and positiveness
(negativeness) for Hermitian subspaces of product spaces of Hilbert spaces,
give their characterizations by induced operators ${\ti T_s}, {\hat T_s}$, and
$T_s$, and finally  establish a close relationship between the boundedness
and boundedness both from below and from above for a Hermitian subspace.

Let $X$ be an inner product space. A subspace $T\in LR(X)$
is said to be Hermitian if $\la y_2, x_1 \ra=\la x_2, y_1 \ra$
for any $(x_1,y_1), (x_2,y_2)$ $\in T$ [3].
We first give some properties of Hermitian subspaces.\medskip

\noindent{\bf Proposition 3.1.} Let $X$ be an inner product space and $T\in LR(X)$
be Hermitian. Then $D(T)\subset {\ol {T(0)}}^\perp$. Further,
$T$ is single-valued if $D(T)$ is dense in $X$.\medskip

\noindent{\bf Proof.} For any given $x\in D(T)$, there exists $y\in X$
such that $(x,y)\in T$. Since $T$ is Hermitian, we have that
$\la z, x \ra=\la 0, y \ra=0$ for all $z\in T(0),$
which yields that $x\in T(0)^\perp={\ol {T(0)}}^\perp$. Hence,
$D(T)\subset {\ol {T(0)}}^\perp$.

Further, suppose that $D(T)$ is dense in $X$. Then it follows from the above assertion
that ${\ol {T(0)}}
\subset D(T)^\perp=\{0\}$, which implies that $T(0)=\{0\}$.
Therefore, $T$ is single-valued. The proof is complete.\medskip

\noindent{\bf Proposition 3.2.} Let $X$ be an inner product space and $T\in LR(X)$.
Then $T$ is Hermitian if and only if so is ${\ol T}$.\medskip

\noindent{\bf Proof.} The sufficiency is obvious. Now, we show the necessity.
Suppose that $T$ is Hermitian.
Fix any $(x_1,y_1), (x_2,y_2)\in {\ol T}$. There exist two sequences
$\{(x_j^{(n)},y_j^{(n)})\}_{n=1}^\infty\subset T$ such that $(x_j^{(n)},y_j^{(n)})
\to (x_j,y_j)$ as $n\to \infty$ for $j=1,2$. Then
\vspace{-0.2cm}$$\la y_1^{(n)}, x_2^{(n)}\ra
=\la x_1^{(n)}, y_2^{(n)}\ra,\;\;n\ge 1.                             \vspace{-0.2cm}$$
Letting $n\to \infty$ in the above relation, we get that
$\la y_1, x_2\ra =\la x_1, y_2\ra.$
Hence, ${\ol T}$ is Hermitian. The proof is complete.\medskip

\noindent{\bf Proposition 3.3.} Let $X$ be a Hilbert space
and $T\in LR(X)$. Then for all $(x_1,y_1),(x_2,y_2)\in T$,
\vspace{-0.2cm}$$\la {\ti T_s}(x_2),  [x_1]\ra
=\la {\hat T_s}(x_2),  x_1\ra=\la y_2,  x_1\ra.                                    \eqno (3.1)\vspace{-0.2cm}$$
Further, if $T$ is Hermitian, then
\vspace{-0.2cm}$$\la {\ti T_s}(x_2),  [x_1]\ra=\la [x_2],  {\ti T_s}(x_1)\ra,\;\;
\la {\hat T_s}(x_2),  x_1\ra=\la x_2,  {\hat T_s}(x_1)\ra,                               \eqno (3.2)\vspace{-0.2cm}$$
and consequently ${\hat T_s}$ is a Hermitian operator in $X$.
\medskip

\noindent{\bf Proof.} Fix any $(x_1,y_1), (x_2,y_2)\in T$.
Then ${\ti T_s}(x_j)=[y_j]\in X/{\ol {T(0)}}$
for $j=1,2$. In addition, there exist $y_{j,0}\in {\ol {T(0)}}$ and
$y_j^\perp\in {\ol {T(0)}}^\perp$ such that $y_j=y_{j,0}+y_j^\perp$ for $j=1,2$.
Then ${\hat T_s}(x_j)=y_j^\perp$ for $j=1,2$.
Noting that $x_1,x_2\in {\ol {T(0)}}^\perp$ by Proposition 3.1,
by (2.1) we have that
\vspace{-0.2cm}$$\la {\ti T_s}(x_2),  [x_1]\ra=\la [y_2],  [x_1]\ra=\la y_2^\perp,  x_1\ra
=\la y_2,  x_1\ra.                                                             \vspace{-0.2cm}$$
In addition, we get that
\vspace{-0.2cm}$$\la {\hat T_s}(x_2),  x_1\ra=\la y_2^\perp,  x_1\ra=\la y_2,  x_1\ra.                   \vspace{-0.2cm}$$
Hence, (3.1) holds.

Similarly, one has that $\la [x_2],{\ti T_s}(x_1)\ra=\la x_2,  y_1\ra$
and $\la x_2,{\hat T_s}(x_1)\ra=\la x_2,  y_1\ra.$
Further, if $T$ is Hermitian, then $\la y_2, x_1 \ra=\la x_2, y_1 \ra$.
Thus, (3.2) holds, and consequently ${\hat T_s}$ is  Hermitian.
 This completes the proof.
\medskip

In Section 2.3, we have got some equivalent conditions
for the boundedness of subspaces in Theorem 2.4.
Now, we shall further give some sufficient and necessary conditions for the boundedness
and estimations of the bound.

Let $X$ be a Hilbert space and $T\in LR(X)$.
We introduce the following three constants related to
${\ti T_s}$, ${\hat T_s}$, and $T_s$, separately:
\vspace{-0.2cm}$${\ti C}(T):
=\sup\{|\la [x],{\ti T_s}(x)\ra|:\;x\in D(T)\;{\rm with}\;\|x\|\le 1\},            \eqno (3.3)\vspace{-0.2cm}$$
where $[x]\in X/{\ol {T(0)}}$,
\vspace{-0.2cm}$${\hat C}(T):
=\sup\{|\la x,{\hat T_s}(x)\ra|:\;x\in D(T)\;{\rm with}\;\|x\|\le 1\},               \eqno (3.4)\vspace{-0.2cm}$$
and further, in the case that $T\in CR(X)$, denote
\vspace{-0.2cm}$$ C(T):
=\sup\{|\la x, T_s(x)\ra|:\;x\in D(T)\;{\rm with}\;\|x\|\le 1\}.                             \eqno (3.5)$$

\noindent{\bf Proposition 3.4.} Let $X$ be a Hilbert space and $T\in LR(X)$.
Then
\vspace{-0.2cm}$${\ti C}(T)={\hat C}(T).                                                \eqno (3.6) \vspace{-0.2cm}$$
Further, if $T$ is closed, then
\vspace{-0.2cm}$${\ti C}(T)={\hat C}(T)=C(T).                                                   \eqno (3.7)$$

\noindent{\bf Proof.}  (3.6) directly follows from (3.1), (3.3), and (3.4).
In the case that $T$ is closed, $T_s={\hat T}_s$ by Theorem 2.2.
Hence, (3.7) holds by (3.5) and (3.6). This completes the proof.
\medskip

\noindent{\bf Theorem 3.1.} Let $X$ be a complex Hilbert space and $T\in LR(X)$
with dense domain $D(T)$ in $X$. Then\vspace{-0.2cm}
\begin{itemize}
\item[{\rm (i)}] the subspace $T$ is bounded if and only if ${\ti C}(T)<\infty$
and $\|T\|\le 2\, {\ti C}(T)$;\vspace{-0.2cm}
\item[{\rm (ii)}] the subspace $T$ is bounded if and only if
${\hat C}(T)<\infty$ and $\|T\|\le 2\, {\hat C}(T)$;\vspace{-0.2cm}
\item[{\rm (iii)}] the subspace $T$ is bounded if and only if
$C(T)<\infty$ and $\|T\|\le 2\, C(T)$ in the case that $T$ is closed.
\end{itemize}

\noindent{\bf Proof.} We first consider Assertion (ii).
By the assumption that $D(T)$ is dense in $X$ and the fact that $D({\hat T_s})=D(T)$,
${\hat T_s}$ is a densely defined operator in $X$.
Hence, we get by [36, Theorem 4.4]
that ${\hat T_s}$ is bounded if and only if ${\hat C}(T)<\infty$
and $\|{\hat T_s}\|\le 2 \,{\hat C}(T)$.
Consequently, Assertion (ii) follows from (2.15).

Assertions (i) and (iii) are directly derived from Proposition 3.4
and the above conclusion. The entire proof is complete.\medskip

\noindent{\bf Remark 3.1.} The results of Theorem 3.1 extend
those of [36, Theorem 4.4] for operators to subspaces.\medskip

Next, we shall recall concepts of boundedness from below (above),
non-negativeness (non-positiveness), and positiveness
(negativeness) for Hermitian subspaces.
\medskip

\noindent{\bf Definition 3.1 {\rm [33, Definition 2.4]}.} Let $X$ be a Hilbert space and $T\in LR(X)$
be Hermitian.
\begin{itemize}\vspace{-0.2cm}
\item [{\rm (1)}] $T$ is said to be bounded from below (above)
if there exists a constant $C\in
{\mathbf R}$ such that
\vspace{-0.2cm}$$\la y, x\ra \geq C\|x\|^2\;(\langle y, x\rangle \le C\|x\|^2),
 \; \,(x,y)\in T,                                                                          \vspace{-0.1cm}$$
while such a constant $C$ is called a lower (upper) bound of $T$.\vspace{-0.2cm}

\item [{\rm (2)}] $T$ is said to be non-negative (non-positive)
if $0$ is a lower (upper) bound of $T$.\vspace{-0.2cm}

\item [{\rm (3)}] $T$ is said to be positive (negative) if
\vspace{-0.2cm}$$\la y, x\ra >0\;(\la y, x\ra <0), \; \; \,(x,y)\in T\;{\rm with}\;x\neq 0.                         $$
\end{itemize}

If $T$ is Hermitian, then $\la {\ti T_s}(x), [x] \ra$ is a real value
for any $x\in D(T)$ by Proposition 3.3.
The following result gives equivalent characterizations for boundedness from below (above),
non-negativeness (non-positiveness), and positiveness
(negativeness) for a Hermitian subspace $T$ by ${\ti T_s}$.
It can be directly derived from (3.1) and Definition 3.1.\medskip

\noindent{\bf Theorem 3.2.} Let $X$ be a Hilbert space and $T\in LR(X)$
be Hermitian. Then\vspace{-0.2cm}
\begin{itemize}
\item [{\rm (i)}] $T$ is bounded from below (above)
if and only if there exists a constant $C\in
{\mathbf R}$ such that
\vspace{-0.2cm}$$\la {\ti T_s}(x), [x]\ra \geq C\|x\|^2\;(\la {\ti T_s}(x),[x] \ra \le C\|x\|^2),
 \; \,x\in D(T); \vspace{-0.2cm}                                                                        \vspace{-0.2cm}$$
\item [{\rm (ii)}] $T$ is non-negative (non-positive) if and only if
\vspace{-0.2cm}$$\la {\ti T_s}(x), [x]\ra \ge 0\;(\la {\ti T_s}(x),[x] \ra \le 0),
 \; \,x\in D(T);  \vspace{-0.2cm}                                                \vspace{-0.2cm}$$
\item [{\rm (iii)}] $T$ is positive (negative) if and only if
\vspace{-0.2cm}$$\la {\ti T_s}(x), [x]\ra >0\;(\la {\ti T_s}(x),[x] \ra <0),
 \; \,x\in D(T)\;{\rm with}\;x\neq 0.                                                  $$
\end{itemize}

The following two results give equivalent characterizations for boundedness from below (above),
non-negativeness (non-positiveness), and positiveness
(negativeness) for a Hermitian subspace $T$ by ${\hat T_s}$ and $T_s$, respectively.
They are direct consequences of Proposition 3.3
and Theorem 2.2.\medskip

\noindent{\bf Theorem 3.3.} Let $X$ be a Hilbert space and $T\in LR(X)$
be Hermitian. Then\vspace{-0.2cm}
\begin{itemize}
\item [{\rm (i)}] $T$ is bounded from below (above)
if and only if so is ${\hat T_s}$ with the same lower (upper) bound;\vspace{-0.2cm}

\item [{\rm (ii)}] $T$ is non-negative (non-positive)
if and only if so is ${\hat T_s}$;\vspace{-0.2cm}

\item [{\rm (iii)}] $T$ is positive (negative) if and only if so is ${\hat T_s}$.
\end{itemize}

\noindent{\bf Theorem 3.4.} Let $X$ be a Hilbert space and $T\in CR(X)$
be Hermitian. Then\vspace{-0.2cm}
\begin{itemize}
\item [{\rm (i)}] $T$ is bounded from below (above)
if and only if so is $T_s$ with the same lower (upper) bound;\vspace{-0.2cm}

\item [{\rm (ii)}] $T$ is non-negative (non-positive)
if and only if so is $T_s$;\vspace{-0.2cm}

\item [{\rm (iii)}] $T$ is positive (negative) if and only if so is $T_s$.
\end{itemize}

\noindent{\bf Remark 3.2.} Assertion (i) of Theorem 3.4 is the same as that
of [33, Lemma 2.4].\medskip

At the end of this subsection, we shall give a close relationship between
boundedness and boundedness both from below and from above for Hermitian subspaces.
\medskip

\noindent{\bf Theorem 3.5.} Let $X$ be a Hilbert space and $T\in LR(X)$
be Hermitian. If $T$ is bounded, then $T$ is bounded both from below
and from above. In addition, the converse conclusion holds in the case that $D(T)$ is dense
in $X$.\medskip

\noindent{\bf Proof.} Since $T$ is Hermitian, ${\hat T_s}$ is a Hermitian operator
by Proposition 3.3. Note that $D({\hat T_s})=D(T)$.

Suppose that $T$ is bounded. Then $\|T\|<+\infty$.
It follows from Theorem 2.4 and (ii) of Lemma 2.6 that
\vspace{-0.2cm}$$|\la{\hat T_s}(x), x\ra|\le \|{\hat T_s}(x)\|\|x\|
=\|T(x)\|\|x\|\le \|T\|\|x\|^2, \;\; x\in D(T),\vspace{-0.2cm}$$
which implies that ${\hat T_s}$ is bounded both from below
and from above. Hence, $T$ is bounded both from below
and from above by (i) of Theorem 3.3.

Now, suppose that $T$ is bounded both from below
and from above and $D(T)$ is dense
in $X$. Again by (i) of Theorem 3.3, ${\hat T_s}$ is bounded both from below
and from above. So there exist real constants $C_1$ and $C_2$
such that
\vspace{-0.2cm}$$C_1\|x\|^2\le \la {\hat T_s}(x),x\ra \le C_2\|x\|^2,\;\; x\in D(T),\vspace{-0.2cm}$$
which yields that
\vspace{-0.2cm}$$|\la {\hat T_s}(x), x\ra|\le \max\{|C_1|,|C_2|\}\|x\|^2,\;\; x\in D(T).\vspace{-0.1cm}$$
Hence, we get that ${\hat C}(T)\le \max\{|C_1|,|C_2|\}<+\infty$.
Consequently, $T$ is bounded by (ii) of Theorem 3.1. This completes
the proof.\medskip

\noindent{\bf Remark 3.3.} In 2013, we gave another definition of boundedness
for a Hermitian subspace $T$ in (3) of Definition 2.4 in [33]; that is,
$T$ was said to be  bounded if it is bounded both from below
and from above. It follows from Theorem 3.5 that this definition
is weaker than that in Definition 2.1 (see [33, Remark 2.4] for a counterexample),
and they are equivalent in the case that
its domain $D(T)$ is dense in $X$.
We shall remark that the definition in Definition 2.1 is more reasonable.
\bigskip

\noindent{\bf 3.2. Relative boundedness for subspaces}\medskip

In this subsection, we shall introduce a concept of relative boundedness for
subspaces and give its equivalent characterizations by their induced operators.
\medskip

\noindent{\bf Definition 3.2 {\rm [9, Definition VII.2.1]}.} Let $X$, $Y$, and $Z$ be normed spaces,
$T\in LR(X,Y)$, and $S\in LR(X,Z)$.\vspace{-0.2cm}
\begin{itemize}
\item [{\rm (1)}] The subspace $S$ is said to be $T$-bounded if $D(T)\subset D(S)$ and
there exists a constant $c\ge 0$ such that
\vspace{-0.2cm}$$\|S(x)\|\le c\,(\|x\|+\|T(x)\|),\;\;\,x\in D(T).                                      \vspace{-0.2cm}$$
\item [{\rm (2)}] If $S$ is $T$-bounded, then the infimum of all numbers $b\ge 0$
for which a constant $a\ge 0$ exists such that
\vspace{-0.2cm}$$\|S(x)\|\le a\, \|x\|+b\,\|T(x)\|,\;\;\,x\in D(T),                                \vspace{-0.2cm}$$
is called the $T$-bound of $S$.
\end{itemize}

\noindent{\bf Remark 3.4.} In 2014, we gave a definition of relative boundedness
for subspaces (see Definition 2.3
in [32]); that is, $S$ was said to be $T$-bounded if $D(T)\subset D(S)$ and
there exists a constant $c\ge 0$ such that for all $(x,y)\in T$ and $(x,z)\in S$,
\vspace{-0.2cm}$$\|z\|\le c(\|x\|+\|y\|).                                                 \vspace{-0.2cm}$$
The above condition implies that $S$ is single-valued whether $T$ is single-valued
 or multi-valued. So this definition is not reasonable in the case that
$S$ is multi-valued. We shall take this opportunity to express our apology
for our carelessness!
Here, the concept is defined by the norms of $T$ and $S$ at points in Definition 3.2
so that the influence of the multi-valued parts of $T$ and $S$ has been removed.
\medskip

\noindent{\bf Theorem 3.6.} Let $X$, $Y$, and $Z$ be normed spaces,
$T\in LR(X,Y)$, and $S\in LR(X,Z)$. Then\vspace{-0.2cm}
\begin{itemize}
\item [{\rm (i)}] the subspace $S$ is  $T$-bounded with $T$-bound $b$
if and only if ${\ti S}_s$ is ${\ti T}_s$-bounded with ${\ti T_s}$-bound $b$;\vspace{-0.2cm}

\item [{\rm (ii)}] in the case that $Y$ and $Z$ are Hilbert spaces,
the subspace $S$ is  $T$-bounded with $T$-bound $b$
if and only if ${\hat S}_s$ is ${\hat T}_s$-bounded with ${\hat T}_s$-bound $b$;\vspace{-0.2cm}

\item [{\rm (iii)}] in the case that $X=Y=Z$ is a Hilbert space and $T,S\in CR(X)$,
the subspace $S$ is  $T$-bounded with $T$-bound $b$
if and only if $S_s$ is $T_s$-bounded with $T_s$-bound $b$.
\end{itemize}

\noindent{\bf Proof.} Assertion (i) and Assertions (ii)-(iii)
directly follow from (2.9) and Theorem 2.4, respectively. The proof is complete.\medskip

\noindent{\bf Remark 3.5.} More recently, we gave stability of self-adjointness
of subspaces  under the assumption of relative perturbation of their induced operator parts,
introduced by Arens [3], with relative bound less than $1$ (see [31, Theorem 4.1]).
By (iii) of Theorem 3.6, the relative perturbation of two closed subspaces $T$ and $S$
in a Hilbert space $X$ is the same as that of their induced operator parts $T_s$ and $S_s$.

\bigskip

\noindent{\bf 4. Closedness and closability of subspaces}
\medskip

In this section, we shall discuss properties and characterizations
for the closedness and closability of subspaces. In particular,
we shall consider relationships of colsedness (closability) of subspace $T$ with
that of its induced operators ${\ti T_s}$ and ${\hat T_s}$.

We first give the following fundamental results:\medskip

\noindent{\bf Lemma 4.1.}
Let $X$ be a topological linear space, $Y$ be a normed space, and $T\in LR(X,Y)$.
Then ${\ol {({\ti T_s})}}=G(Q_T){\ol T}.$
\medskip

\noindent{\bf Proof.} The proof is similar to that of [9, Proposition II.5.2].
So its details are omitted.\medskip

\noindent{\bf Lemma 4.2.}
Let $X$ and $Y$ be topological linear spaces and $T\in LR(X,Y)$.
Then the following statements are equivalent:
\begin{itemize}\vspace{-0.2cm}
\item [{\rm (i)}] $T$ is closable;\vspace{-0.2cm}
\item [{\rm (ii)}] $T(0)={\ol T}(0)$;\vspace{-0.2cm}
\item [{\rm (iii)}] ${\ti T_s}$ is closable and $T(0)$ is closed
in the case that $Y$ is a normed space.
\end{itemize}
\medskip

\noindent{\bf Proof.} We first show that the statements (i) and (ii) are equivalent.
Suppose that (i) holds. Then $T(x)={\ol T}(x)$ for every $x\in D(T)$.
Thus, (ii) holds. Conversely, we suppose that (ii) holds.
Fix any $x\in D(T)$ and any $y\in T(x)$. Then $y\in {\ol T}(x)$.
So by Lemma 2.1 we have that
\vspace{-0.2cm}
$${\ol T}(x)=\{y\}+{\ol T}(0)=\{y\}+ T(0)=T(x).\vspace{-0.1cm}$$
Hence, $T$ is closable, and then (i) holds.

Now, suppose that $Y$ is a normed space. With a similar argument
to that used in the proof of [9, Proposition II.5.7],
one can show that the statements (i) and (iii) are equivalent.
The proof is complete.\medskip

\noindent{\bf Remark 4.1.} Cross gave the same results as those in Lemmas 4.1 and 4.2
in the case that $X$ and $Y$ are normed spaces (see [9,
Propositions II.5.2, II.5.5, and II.5.7]).

\medskip

\noindent{\bf Lemma 4.3.} Let $X$, $Y$, and $Z$ be topological linear spaces.
Assume that $T\in LR(X,Y)$, $S\in LR(X,Z)$,
$U:\;Y\to Z$ is a linear and homeomorphic operator, and they satisfy
that $D(T)=D(S)$ and
\vspace{-0.2cm}$$S=G(U)\,T.                                                                   \eqno (4.1) \vspace{-0.2cm}$$
Then $T$ is closed (closable) if and only if $S$ is closed (closable).\medskip

\noindent{\bf Proof.} It follows from (4.1) that ${\ol S}=G(U)\,{\bar T}$.
So it suffices to show that the assertion about the closedness holds.
In addition, we have that $T=\lf(G(U)\rf)^{-1}S=G(U^{-1})\,S$ by (4.1).
Thus, it is enough for us to show that the sufficiency of the assertion
about the closedness holds.

Suppose that $S$ is closed. Let $\{(x_n,y_n)\}_{n=1}^\infty\subset T$ be any convergent
sequence with $x_n\to x$ and $y_n\to y$ as $n\to \infty$. Let $z_n=Uy_n$ for $n\ge 1$
and $z=Uy$.
Then $(x_n,z_n)\in S$ for $n\ge 1$. By the assumption that $U$ is homeomorphic, one has
that $z_n\to z$ as $n\to \infty.$
Hence, $(x,z)\in S$, and consequently $x\in D(T)$ by the assumption that $D(S)=D(T)$
and $(x,y)\in T$. Therefore, $T$ is closed. This completes the proof.
\medskip

Let $T\in LR(X,Y)$. In the case that $X$ is a topological linear space
and $Y$ is a normed space,
Lemmas 2.2 and 4.2 give out equivalent characterizations for closedness and
closability of the subspace $T$ by ${\ti T_s}$
and $T(0)$, respectively.
The following result gives another equivalent characterization
of closedness and closability  of the subspace $T$ by ${\hat T_s}$ and $T(0)$, and the equivalence
between closedness (closability) of ${\ti T_s}$ and that of ${\hat T_s}$.\medskip

\noindent{\bf Theorem 4.1.} Let $X$ be a topological linear space, $Y$ be a Hilbert space,
and $T\in LR(X,Y)$. Then\vspace{-0.2cm}
\begin{itemize}
\item [{\rm (i)}] ${\ti T_s}$ is closed (closable) if and only if so is ${\hat T_s}$;\vspace{-0.2cm}
\item [{\rm (ii)}] $T$ is closed (closable) if and only if ${\hat T_s}$ is closed (closable)
and $T(0)$ is closed.
\end{itemize}

\medskip

\noindent{\bf Proof.} It can be easily verified that Assertion (i)
holds by Lemmas 2.4 and 4.3 and Theorem 2.1. So Assertion (ii) holds by Lemmas 2.2
and 4.2. The proof is compete.\medskip

The following two results give other two equivalent characterizations for
closedness of a subspace in terms of its graph norm and graph inner product.

Let $X$ and $Y$ be normed spaces and $T\in LR(X,Y)$. Define
\vspace{-0.2cm}$$\|x\|_T:=\|x\|+\|T(x)\|,\;\;x\in D(T).                                  \eqno (4.2) \vspace{-0.2cm}$$
We call $\|\cdot\|_T$ defined by (4.2) the graph norm for the subspace $T$.\medskip

We shall first recall the following result for operators in Banach spaces:\medskip

\noindent{\bf Lemma 4.4 {\rm [10, Theorem 1.3.1]}}. Let $X$ and $Y$ be Banach spaces
and $T$ be an operator from $X$ to $Y$. Then
$T$ is closed if and only if $(D(T),\|\cdot\|_T)$ is a Banach space.
\medskip

\noindent{\bf Theorem 4.2.} Let $X$ and $Y$ be Banach spaces and $T\in LR(X,Y)$. Then
$T$ is closed if and only if $(D(T),\|\cdot\|_T)$ is a Banach space and $T(0)$ is closed.
\medskip

 \noindent{\bf Proof.} Note that  $D(T)=D({\ti T_s})$ and
 $\|\cdot\|_T=\|\cdot\|_{\ti T_s}$ by (2.9). Hence, $(D(T),\|\cdot\|_T)$ is a Banach space
 if and only if so is $(D({\ti T_s}),\|\cdot\|_{\ti T_s})$. In addition,
 $Y/T(0)$ is a Banach space if $T(0)$ is closed, and ${\ti T_s}$ is single-valued.
 One can easily show that the assertion in Theorem 4.2 holds by Lemmas 2.2 and 4.4.
 The proof is complete.\medskip

 \noindent{\bf Theorem 4.3.} Let $X$ and $Y$ be Hilbert spaces and $T\in LR(X,Y)$. Then
$T$ is closed if and only if $(D(T),\la \cdot,\cdot \ra _T)$ is a Hilbert space
and $T(0)$ is closed,
where
\vspace{-0.2cm}$$\la x_1, x_2\ra _T: = \la x_1, x_2\ra + \la T(x_1), T(x_2)\ra ,\;\;x_1,x_2\in D(T),   \eqno (4.3) \vspace{-0.2cm}$$
while
\vspace{-0.2cm}$$\la T(x_1), T(x_2)\ra:=\la [y_1], [y_2]\ra,\;\; y_1\in T(x_1),\; y_2\in T(x_2),        \eqno (4.4) \vspace{-0.2cm}$$
and $\la [y_1], [y_2]\ra$ is defined by (2.1).\medskip

We shall remark that the inner product in (4.4) is well defined since
$\la [y_1], [y_2]\ra= \la [y_1'], [y_2']\ra$ for any $y_j,y_j'\in T(x_j)$, $j=1,2$.
We call $\la \cdot, \cdot\ra_T$ defined by (4.3) and (4.4) the graph inner
product for the subspace $T$.
\medskip

\noindent{\bf Proof.} It is evident that $(D(T),\la \cdot,\cdot \ra _T)$
is an inner product space, and the induced norm by
$\la \cdot, \cdot\ra _T$
is equivalent to the graph norm $\|\cdot\|_T$ for the subspace $T$.
In addition, $(D(T),\la \cdot,\cdot \ra _T)$ is a Hilbert space if and only if
$(D(T),\|\cdot\|_T)$ is a Banach space. Hence, this theorem
follows from Theorem 4.2. This completes the proof.\medskip

The following result gives a sufficient condition for closability of subspaces.
\medskip

\noindent{\bf Theorem 4.4.} Let $X$ and $Y$ be normed spaces and $T\in LR(X,Y)$.
If $T$ is bounded and $T(0)$ is closed, then $T$ is closable.
Further, if $Y$ is complete, then
\vspace{-0.2cm}$$D({\ol T})={\ol {D(T)}},                                              \eqno (4.5)\vspace{-0.2cm}$$
 and ${\ol T}$ is the bounded extension of $T$
onto ${\ol {D(T)}}$ with bound $\|T\|$.\medskip

\noindent{\bf Proof.} Since $T$ is bounded, $T$ is continuous by [9, (a) of Proposition II.3.2].
So $T$ is closable by [9, Proposition II.5.7].

Further, suppose that $Y$ is complete. Now, we show that (4.5) holds.

We shall first consider the special case that $T$ is single-valued.
It is evident that $D({\ol T})\subset {\ol {D(T)}}$ by the definition of ${\ol T}$.
So it is only needed to show that $D({\ol T})\supset {\ol {D(T)}}$.
For any given $x\in {\ol {D(T)}}$, there exists a sequence
$\{x_n\}_{n=1}^\infty\subset D(T)$ such that $x_n\to x$
as $n\to \infty$. Then $\{T(x_n)\}_{n=1}^\infty$ is a Cauchy sequence
in $Y$ since $T$ is a bounded operator. Hence, there exists $y\in Y$ such that
$T(x_n)\to y$ as $n\to \infty$. Thus, $(x, y)\in {\ol {T}}$,
which implies that $x\in D({\ol T})$.
And consequently, $D({\ol T})\supset {\ol {D(T)}}$. Therefore, (4.5) holds in this case.

Next, we consider the general case that $T$ is multi-valued.
Since $T$ is bounded and $T(0)$ is closed,
${\ti T_s}$ is a bounded operator from $X$ to $Y/T(0)$ by (2.9)
and $Y/T(0)$ is complete. Hence, $D({\ol {({\ti T_s})}})={\ol {D({\ti T_s})}}$
by the above discussion. In addition, $D({\ol {({\ti T_s})}})=D({\ol T})$
by Lemma 4.1 and $D({\ti T_s})=D(T)$. Therefore, (4.5) holds in the general case.

It follows from [9, Proposition II.5.8] that $\|{\ol T}\|= \|T\|$ if $T$ is closable.
Hence, the last statement in Theorem 4.4 holds. The whole proof is complete.\medskip

At the end of this section, we give another equivalent characterization for
closability of a subspace.\medskip

\noindent{\bf Theorem 4.5.} Let $X$ and $Y$ be topological linear spaces and $T\in LR(X,Y)$.
Then $T$ is closable if and only if $T$ has a closed extension; that is, there
exists $S\in CR(X, Y)$ such that $D(T)\subset D(S)$
and  $S(x)=T(x)$ for any $x\in D(T)$.\medskip

\noindent{\bf Proof.} The necessity is obvious. Now, we show the sufficiency.
Suppose that $T$ has a closed extension $S\in CR(X, Y)$. Then
$T\subset S$ and $S(x)=T(x)$ for any $x\in D(T)$. So ${\ol T}\subset S$,
and then
\vspace{-0.2cm}$${\ol T}(x)\subset S(x)=T(x), \;\; x\in D(T),\vspace{-0.1cm}$$
which, together with the fact that $T(x)\subset {\ol T}(x)$ for every $x\in D(T)$,
implies that ${\ol T}(x)=T(x)$ for every $x\in D(T)$. Hence, $T$ is closable.
This completes the proof.\medskip

\noindent{\bf Remark 4.2.} The results of Theorems 4.2-4.5 extend those for operators
in Hilbert spaces  and Banach spaces to subspaces in  Banach spaces or normed spaces or topological linear spaces,
respectively (see [10, Theorem 1.3.1] and [36, Theorems 5.1, 5.2, and 5.4]).\medskip

\bigskip

\noindent{\bf 5. Closed graph theorem}\medskip

It is well known that the closed graph theorem for operators
plays an important role in the study of linear operators
in Hilbert and Banach spaces (cf., [14, 36]).
So is it for subspaces in the study of subspaces.
In this section, we shall focus our attention on this subject.

We first recall the following result for operators in Banach spaces:\medskip

\noindent{\bf Lemma 5.1 {\rm [14, Chap. III, Theorem 5.20]}.}
Let $X$ and $Y$ be Banach spaces. If $T$ is a closed
operator from $X$ to $Y$ and its domain $D(T)$ is closed, then $T$ is bounded.
\medskip

Note that in [14, Chap. III, Theorem 5.20], it is required that $D(T)=X$. It is evident
that the above result still holds since $D(T)$ can be regarded as a Banach space.
\medskip

\noindent{\bf Theorem 5.1 {\rm (Closed graph theorem for subspaces in Banach spaces).}}
Let $X$ and $Y$ be Banach spaces and $T\in LR(X,Y)$. Then the following statements are equivalent:
\begin{itemize}\vspace{-0.2cm}
\item [{\rm (i)}] $T$ and $D(T)$ are closed;\vspace{-0.2cm}
\item [{\rm (ii)}] $T$ is bounded, and $D(T)$ and $T(0)$ are closed;\vspace{-0.2cm}
\item [{\rm (iii)}] $T$ is bounded and closed.
\end{itemize}

\noindent{\bf Proof.} ``(i) $\Rightarrow$ (ii)". Suppose that (i) holds.
Then ${\ti T_s}$ and $T(0)$ are closed by Lemma 2.2. Noting that
$Y/\,{\ol {T(0)}}=Y/T(0)$ is a Banach space and $D({\ti T_s})=D(T)$,
we get that ${\ti T_s}$ is bounded by Lemma 5.1. Hence, $T$ is bounded by (2.9), and
then (ii) holds.

``(ii) $\Rightarrow$ (iii)". Suppose that (ii) holds. It suffices to show that
${\ti T_s}$ is closed by Lemma 2.2. Fix any convergent sequence
$\{(x_n,[y_n])\}_{n=1}^\infty\subset {\ti T_s}$
with $\{(x_n,y_n)\}_{n=1}^\infty\subset T$, and $x_n\to x$ and $[y_n]\to [y]$
as $n\to \infty$. Then $x\in D(T)$ since $D(T)$ is closed. For any given
$y'\in T(x)$,  we have that $(x,[y'])\in {\ti T_s}$.
Since $T$ is bounded, ${\ti T_s}$ is bounded again
by (2.9). Then
\vspace{-0.2cm}$$\|[y_n]-[y']\|=\|{\ti T_s}(x_n)-{\ti T_s}(x)\|\le \|{\ti T_s}\|\,\|x_n-x\|,\;n\ge 1,\vspace{-0.2cm}$$
which implies that $[y_n]\to [y']$ as $n\to \infty$. Hence, $[y]=[y']$, and consequently
$(x,[y])\in {\ti T_s}$. Therefore, ${\ti T_s}$ is closed and (iii) holds.

``(iii) $\Rightarrow$ (i)" Suppose that (iii) holds. Then ${\ti T_s}$ is bounded, and
${\ti T_s}$ and $T(0)$ are closed
by Lemma 2.2 and (2.9). Obviously, it is only needed to show that
$D(T)$ is closed. Fix any convergent sequence
$\{x_n\}_{n=1}^\infty\subset D(T)$
with $x_n\to x$ as $n\to \infty$.
For any given $y_n\in T(x_n)$ for each $n\ge 1$, we have that $(x_n,[y_n])\in {\ti T_s}$
. Since ${\ti T_s}$ is bounded, we get that
\vspace{-0.2cm}$$\|[y_n]-[y_m]\|=\|{\ti T_s}(x_n)-{\ti T_s}(x_m)\|\le \|{\ti T_s}\|\,\|x_n-x_m\|,\;n,m\ge 1.\vspace{-0.2cm}$$
Thus, $\{[y_n]\}_{n=1}^\infty$ is a Cauchy sequence in $Y/\,{\ol {T(0)}}$,
and then is convergent by the fact that $Y/\,{\ol {T(0)}}=Y/T(0)$ is a Banach space.
Denote $[y]=\lim _{n\to \infty}[y_n]$. Then $(x,[y])\in {\ti T_s}$ because
${\ti T_s}$ is closed. Hence, $x\in D({\ti T_s})=D(T)$. Consequently,
$D(T)$ is closed.

The entire proof is complete.
\medskip

\noindent{\bf Remark 5.1.} The result of Theorem 5.1 extends the closed graph theorem
for operators in Hilbert spaces (see [36, Theorem 5.6]) to subspaces in Banach
spaces.
\medskip

The following three results are direct consequences of Theorem 5.1.\medskip

\noindent{\bf Corollary 5.1.}
Let $X$ and $Y$ be Banach spaces and $T\in CR(X,Y)$. Then
$T$ is bounded if and only if $D(T)$ is closed.
\medskip

\noindent{\bf Remark 5.2.} The result of Corollary 5.1 is the same as that
of [9, (a) of Theorem III.4.2 ]. Recently,
 \'Alvarez generalized the sufficiency result
of Corollary 5.1 for the Banach space $Y$ to the paracomplete space $Y$
[1, Proposition 8].\medskip

\noindent{\bf Corollary 5.2.}
Let $X$ and $Y$ be Banach spaces and $T\in LR(X,Y)$ with closed domain $D(T)$. Then
$T$ is closed if and only if $T$ is bounded and $T(0)$ is closed.
\medskip

\noindent{\bf Corollary 5.3.}
Let $X$ and $Y$ be Banach spaces and $T\in BR(X,Y)$. Then
$T$ is closed if and only if $D(T)$ and $T(0)$ are closed.
\medskip

\noindent{\bf Theorem 5.2.} Let $X$ and $Y$ be Banach spaces and
$T\in LR(X,Y)$ with closed range $R(T)$. Then
$T$ is closed if and only if $T^{-1}$ is bounded and $N(T)$ is closed.\medskip

\noindent{\bf Proof.} Note that
\vspace{-0.2cm}$$D(T^{-1})=R(T),\;\;T^{-1}(0)=N(T).\vspace{-0.1cm}$$
So $D(T^{-1})$ is closed by the assumption that $R(T)$ is closed.
Hence, $T^{-1}$ is closed if and only if $T^{-1}$ is bounded
 and $N(T)$ is closed by Corollary 5.2.
Consequently,
the assertion holds by the fact that $T$ is closed if and only if  so is $T^{-1}$.
This completes the proof.

\medskip

\noindent{\bf Remark 5.3.} The result of Theorem 5.2 extends
the open mapping theorem for operators in Hilbert spaces
(see [36, Theorem 5.8]) to subspaces in Banach spaces.
In the operator case, it is required that the operator $T$
is injective. Now, this requirement has been relaxed as
that $N(T)$ is closed.\medskip

\bigskip

\noindent{\bf 6. Stability of closedness and closability for subspaces under perturbation}\medskip

In this section, we shall first give two sufficient conditions
for relative boundedness of a subspace, and then study stability of closedness
and closability for subspaces under bounded and relatively bounded perturbations.\medskip

\noindent{\bf Theorem 6.1.} Let $X$, $Y$, and $Z$ be normed spaces,
$T\in LR(X,Y)$, and $S\in LR(X,Z)$ with $D(T)\subset D(S)$. If $S$
is bounded, then $S$ is $T$-bounded with $T$-bound $0$.\medskip

\noindent{\bf Proof.} Since $S$ is bounded, we have that
$\|S\|<+\infty$ and
\vspace{-0.2cm}$$\|S(x)\|\le \|S\|\|x\|,\;\; x\in D(S),\vspace{-0.2cm}$$
which implies that for any number $b>0$,
\vspace{-0.2cm}$$\|S(x)\|\le \|S\|\|x\|+b \|T(x)\|,\;\; x\in D(T).\vspace{-0.2cm}$$
Hence, $S$ is $T$-bounded with $T$-bound $0$. This completes the proof.
\medskip

\noindent{\bf Theorem 6.2.} Let $X$, $Y$, and $Z$ be Banach spaces,
$T\in LR(X,Y)$, and $S\in LR(X,Z)$ with $D(T)\subset D(S)$.
If $T$ is closed and $S$ is closable, then $S$ is $T$-bounded.
\medskip

\noindent{\bf Proof.} Since $T$ is closed, $T(0)$ is closed
and $(D(T),\|\cdot\|_T)$ is a Banach space
by Theorem 4.2, where $\|\cdot\|_T$ is defined by (4.2).

Define the subspace $S_0\in LR((D(T),\|\cdot\|_T), Z)$
by $S_0(x)=S(x)$ for $x\in D(T)$.
In order to show that
$S$ is $T$-bounded, it suffices to show that $S_0$ is bounded.
By Corollary 5.2, it is enough to show that $S_0$ is closed.
Fix any sequence $\{(x_n,z_n)\}_{n=1}^\infty \subset S_0$ with $x_n\to x$
in norm $\|\cdot\|_T$ and $z_n\to z$ in the norm of $Z$
as $n\to \infty$. Take any $y_n\in T(x_n)$ for each $n\ge 1$
and $y\in T(x)$. Noting that
\vspace{-0.2cm}$$\|x_n-x\|_T=\|x_n-x\|+\|T(x_n)-T(x)\|=\|x_n-x\|+\|[y_n]-[y]\|,\;n\ge 1,\vspace{-0.2cm}$$
we get that $x_n\to x$ in the norm of $X$ and $[y_n]\to [y]$
in the norm of $Y/T(0)$ as $n\to \infty$. In addition,
$S(x_n)=S_0(x_n)$, which implies that $(x_n,z_n)\in S$
 for each $n\ge 1$.
By the assumption that $S$ is closable and $x\in D(T)\subset D(S)$,
we have that $z\in {\bar S}(x)=S(x)=S_0(x)$. Thus,
$S_0$ is closed. This completes the proof.
\medskip

\noindent{\bf Remark 6.1.} [6, Lemma 4.3] gave the same result as that of Theorem 6.2
in the case that $Y=Z$.\medskip

\noindent{\bf Theorem 6.3.} Let $X$ be a normed space, $Y$ be a Banach space,
and $T, S\in LR(X,Y)$ satisfy that $D(T)\subset D(S)$
and $S(0)\subset T(0)$. If $S$ is $T$-bounded with $T$-bound less than $1$,
then $T+S$ is closed (closable) if and only if $T$ is closed (closable).
Moreover, $D({\ol {T+S}})=D({\ol T})$.\medskip

\noindent{\bf Proof.} Since $S$ is $T$-bounded with $T$-bound less than $1$,
there exist constants $a>0$ and $0\le b<1$ such that
\vspace{-0.2cm}$$\|S(x)\|\le a\|x\|+b\|T(x)\|,\;\;x\in D(T).                        \eqno(6.1)\vspace{-0.2cm}$$
Note that $D(T+S)=D(T)$ by the assumption that $D(T)\subset D(S)$.
Then, by (iii) of Lemma 2.5 and (6.1) one has that
\vspace{-0.2cm}$$\|(T+S)(x)\|
\le \|T(x)\|+\|S(x)\|\le a\|x\|+(b+1)\|T(x)\|,\;\;x\in D(T).                     \eqno(6.2)\vspace{-0.2cm}$$
In addition, it follows from the assumption that $S(0)\subset T(0)$ that
\vspace{-0.2cm}$$(T+S)(0)=T(0)+S(0)=T(0)\supset S(0).                                   \eqno(6.3) \vspace{-0.2cm}$$
So, $Y/\,{\ol {(T+S)(0)}}=Y/\,{\ol {T(0)}}$, and by Theorem 2.3 and (6.1) we get that
\vspace{-0.2cm}$$\|(T+S)(x)\|\ge \|T(x)\|-\|S(x)\|\ge (1-b)\|T(x)\|-a\|x\|,\;\;x\in D(T),   \vspace{-0.2cm}$$
which implies that
\vspace{-0.2cm}$$\|T(x)\|
\le \frac{1}{1-b}\left(a\|x\|+\|(T+S)(x)\|\right),\;\;x\in D(T).                  \eqno(6.4)\vspace{-0.2cm}$$

We shall show this theorem by three steps.

{\bf Step 1.} Show the assertion about the closability.

``$\Leftarrow$" Suppose that $T$ is closable. Then ${\ti T_s}$ is closable,
$T(0)$ is closed, and
\vspace{-0.2cm}$$T(0)={\ol T}(0)                                                                       \eqno(6.5)    \vspace{-0.2cm}$$
 by Lemma 4.2. In order to show that $T+S$ is closable, again by Lemma 4.2
it suffices  to show that
 \vspace{-0.2cm}$$(T+S)(0)={\ol {(T+S)}}(0).                                                            \eqno(6.6)     \vspace{-0.2cm}$$
By Lemma 2.3 we have that $(T+S)(0)\subset {\ol {(T+S)(0)}}\subset {\ol {(T+S)}}(0)$.
So it is only needed to show
that \vspace{-0.2cm}$${\ol {(T+S)}}(0)\subset (T+S)(0)=T(0),                                                \eqno(6.7)  \vspace{-0.2cm}$$
where (6.3) has been used. Again by (6.3) one has that
\vspace{-0.2cm}$${\ol {(T+S)(0)}}={\ol {T(0)}}=T(0),\vspace{-0.2cm}$$ which, together with (6.5),
yields that \vspace{-0.2cm}$$Y/\,{\ol {T(0)}}=Y/T(0)=Y/\,{\ol {(T+S)(0)}}=Y/\,{\ol T}(0).\vspace{-0.2cm}$$
Hence,
\vspace{-0.2cm}$$Q_T=Q_{T+S}=Q_{\ol T}.                                                               \eqno(6.8)  $$

Fix any $w\in {\ol {(T+S)}}(0)$. Then there exists
$\{(x_n,w_n)\}_{n=1}^\infty\subset T+S$ such that $x_n\to 0$ and
$w_n\to w$ as $n\to \infty$. And there exist $y_n\in T(x_n)$ and
 $z_n\in S(x_n)$ such that $w_n=y_n+z_n$ for each
$n\ge 1$. For clarity, by $[u]_T$ denote the elements
of $Y/\,{\ol {T(0)}}$. Then $[u]_T=[u]_{T+S}=[u]_{{\ol T}}$
 for every $u\in Y$ by (6.8). Noting that ${\ol {S(0)}}\subset T(0)$,
 we have that
\vspace{-0.2cm}$$\|[u]_S\|=d(u, {\ol {S(0)}})\ge d(u, T(0))=\|[u]_T\|,\;\;u\in Y.        \eqno(6.9)\vspace{-0.2cm}$$
It follows from (6.4) that
\vspace{-0.2cm}$$\|T(x_n-x_m)\|\le
\frac{1}{1-b}\left(a\|x_n-x_m\|+\|(T+S)(x_n-x_m)\|\right), \vspace{-0.2cm}$$
which yields that
\vspace{-0.2cm}$$\|[y_n]_T-[y_m]_T\|\le
\frac{1}{1-b}\left(a\|x_n-x_m\|+\|[w_n]_T-[w_m]_T\|\right), \;\;\;n>m\ge 1.           \eqno(6.10)\vspace{-0.2cm}$$
In addition, by (i) of Lemma 2.5 we have that
\vspace{-0.2cm}$$\|[w_n]_T-[w_m]_T\|=d(w_n-w_m,
T(0))\le \|w_n-w_m\|.                                                                      \eqno(6.11) \vspace{-0.2cm}$$
It follows from (6.10) and (6.11) that $\{[y_n]_T\}$ is a Cauchy sequence in $Y/T(0)$,
and so there exists $[y]_T\in
Y/T(0)$ such that $[y_n]_T\to [y]_T$ as $n\to \infty$ by the fact that $Y/T(0)$
is complete. Then $(0, [y]_T)\in {\ol {\ti T_s}}$ by the fact that
$(x_n,[y_n]_T)\in {\ti T_s}$ for each $n\ge 1$.
Further, by Lemma 4.1 and (6.8) we have that
\vspace{-0.2cm}$${\ol {\ti T_s}}
=G(Q_T) {\ol T}=G(Q_{\ol T}) {\ol T}={\widetilde {({\ol T})_s}}.                     \eqno(6.12)\vspace{-0.2cm}$$
Hence, $[y]_T=[0]_T.$ In addition, it follows from (6.1) that
\vspace{-0.2cm}$$\|[z_n]_S\|=\|S(x_n)\|\le a \|x_n\|+b\|T(x_n)\|
=a \|x_n\|+b\|[y_n]_T\|,\;n\ge 1,  \vspace{-0.2cm}$$
which implies that $[z_n]_S\to [0]_S$ as $n\to \infty$.
By (6.9) we get that
$\|[z_n]_T\|\le \|[z_n]_S\|$ for each $n\ge 1$. Hence, $[z_n]_T\to [0]_T$
as $n\to \infty$, and consequently
$[w_n]_T=[y_n]_T+[z_n]_T\to [0]_T$ as $n\to \infty$. So $[w]_T=[0]_T$,
which implies that $w\in T(0)$. Therefore,
(6.7) holds, and then (6.6) holds. Consequently, $T+S$ is closable.

``$\Rightarrow$" Suppose that $T+S$ is closable. It follows from Proposition 2.1 that
\vspace{-0.2cm}$$T=(T+S)-S,                                                                      \eqno(6.13) \vspace{-0.2cm}$$
where the assumption that $S(0)\subset T(0)$ has been used. With a similar argument to that used in the
above discussion, one can show that $T$ is closable.

{\bf Step 2.} Show that $D({\ol {T+S}})=D({\ol T})$.

Suppose that $T$ (or $T+S$) is closable. Then $T+S$ (or $T$) is closable,
and $T(0)=(T+S)(0)$ is closed. Fix any $x\in
D({\ol {T+S}})$. There exists a convergent sequence
$\{(x_n,w_n)\}_{n=1}^\infty\subset T+S$ such that $x_n\to x$ and
$w_n\to w$ as $n\to \infty$. It is evident that $(x,w)\in {\ol {T+S}}$
and $x_n\in D(T+S)=D(T)$ for each $n\ge 1$. In
addition, there exist $y_n\in T(x_n)$ and $z_n\in S(x_n)$
such that $w_n=y_n+z_n$ for each $n\ge 1$. It follows from
(6.10) and (6.11) that $\{[y_n]_T\}_{n=1}^\infty$ is convergent to
some $[y]_T$. Hence, $(x, [y]_T)\in {\ol {\ti
T_s}}={\widetilde {({\ol T})_s}}$ by (6.12) and the fact
that $(x_n, [y_n]_T)\in {\ti T_s}$. This implies that $x\in
D({\widetilde {({\ol T})_s}})=D({\ol T})$, and consequently
$D({\ol {T+S}})\subset D({\ol T})$.

With a similar argument to that used in the above and using (6.13),
one can show that $D({\ol {T+S}})\supset D({\ol
T})$. Therefore, $D({\ol {T+S}})=D({\ol T})$.

{\bf Step 3.} Show the assertion about the closedness.

 ``$\Leftarrow$" Suppose that $T$ is closed. It is evident that $T$
 is closable. Hence, $T+S$ is closable and $D({\ol {T+S}})=D({\ol T})$ by
 the assertions shown in Steps 1 and 2.
 Fix any convergent sequence
 $\{(x_n,w_n)\}_{n=1}^\infty\subset T+S$
with $x_n\to x$ and $w_n\to w$ as $n\to \infty$. Then $(x,y)\in {\ol {T+S}}$.
Noting that $D({\ol {T+S}})=D({\ol T})=D(T)
=D(T+S)$, we get that $x\in D(T+S)$ and $y\in {\ol {(T+S)}}(x)=(T+S)(x)$, and consequently
$(x,y)\in T+S$. Therefore, $T+S$ is closed.

``$\Rightarrow$" The proof for the necessity is similar to that in Step 1, and thus is omitted.

The entire proof is complete.\medskip

\noindent{\bf Remark 6.2.} The sufficiency for the closedness in
Theorem 6.4 was given in [6, Lemma 2.3]. Now, Theorem 6.3 shows that
the condition is not only sufficient
but also necessary.

\medskip

The following result is a direct consequence of Theorems 6.1 and 6.3.\medskip

\noindent{\bf Corollary 6.1.} Let $X$ be a normed space, $Y$ be a Banach space,
and $T, S\in LR(X,Y)$ satisfy that $D(T)\subset D(S)$ and $S(0)\subset T(0)$.
 If $S$ is bounded,
then $T+S$ is closed (closable) if and only if $T$ is closed (closable).
Moreover, $D({\ol {T+S}})=D({\ol T})$.\medskip

 \noindent{\bf Remark 6.3.} The results of Theorems 6.1-6.3
 extend those for operators  in Hilbert spaces
 (see [36, Proposition in Page 93, and Theorems 5.5 and 5.9])
 to subspaces in Banach or normed spaces.
 \medskip

If the assumption on the domains of $T$ and $S$
is strengthened and the
assumption on $Y$ is weakened, then we get the following result:\medskip

\noindent{\bf Theorem 6.4.} Let $X$ and $Y$ be normed spaces
and $T, S\in LR(X,Y)$ satisfy that ${\ol {D(T)}}\subset D(S)$
and $S(0)\subset T(0)$. If $S$ is bounded,
then $T+S$ is closed (closable) if and only if $T$ is closed (closable).
Moreover, $D({\ol {T+S}})=D({\ol T})$.\medskip

\noindent{\bf Proof.} Since the proof is similar to that of Theorem 6.3,
we omit its details. The proof is complete.\medskip

\noindent{\bf Remark 6.4.} The sufficiency for the closedness in
Theorem 6.4 was given in the special case that $S$ is single-valued
in [9, Exercise II.5.16]. The sufficiency for the closedness and closability
in Theorem 6.4 was given in the special case that $X$ and $Y$ are complete
in [2, Proposition 3.1].
Here, Theorem 6.4 shows that the condition is not only sufficient
but also necessary in a more general case.
 \medskip

\bigskip \noindent{\bf \large References}
\def\hang{\hangindent\parindent}
\def\textindent#1{\indent\llap{#1\enspace}\ignorespaces}
\def\re{\par\hang\textindent}
\noindent \vskip 3mm

\re{[1]} T. \'Alvarez, Small perturbation of normally solvable relations,
Publ. Math. Debrecen 80(1-2)(2012) 155--68.

\re{[2]} T. \'Alvarez, A. Ammar, A. Jeribi, On the essential spectra of some matrix of
linear relations, Math. Meth. Appl. Sci. 37(2014) 620--644.

\re{[3]} R. Arens, Operational calculus of linear relations,
Pac. J. Math. 11(1961) 9--23.

\re{[4]} T. Ya. Azizov, J. Behrndt, P. Jonas, C. Trunk,
Compact and finite rank perturbations of linear
relations in Hilbert spaces,
Integr. Equ. Oper. Theory 63(2009) 151--163.

\re{[5]} T. Ya. Azizov, J. Behrndt, P. Jonas, C. Trunk, Spectral points of definite type
and type $\pi$ for linear operators and relations in Krein spaces,
J. Lond. Math. Soc. 83(2011) 768--788.

\re{[6]} E. Chafai, M. Mnif, Perturbation of normally solvable linear relations
in paracomplete spaces, Linear Algebra Appl. 439(2013) 1875--1885.

\re{[7]}  E. A. Coddington, Extension theory of formally normal
and symmetric subspaces, Mem. Am. Math. Soc. 134(1973).

\re{[8]} A. Coddington, A. Dijksma, Adjoint subspaces in Banach spaces,
with applications to ordinary differential subspaces, Ann. Mat. Pura.
Appl. 118(1978) 1--118.

\re{[9]} R. Cross, Multivalued Linear Operators,
Monographs and Textbooks in Pure and Applied Mathematics, vol. 213, Marcel Dekker, 1998.

\re{[10]} J. Derezi\'nski, Unbounded Linear Operators, Lecture notes, 2013.

\re{[11]} A. Dijksma, H. S. V. de Snoo,
Eigenfunction extensions associated with pairs of ordinary differential expressions,
J. Differ. Equations 60(1985) 21--56.

\re{[12]} S. Hassi, H. de Snoo, One-dimensional graph
perturbations of self-adjoint relations, Ann. Aca. Sci. Fenn.
Math. 20(1997) 123--164.

\re{[13]} S. Hassi, H. de Snoo, F. H. Szafraniec,
Componentwise and cartesian decompositions of linear relations,
Dissertationes Math. 465(2009) (59 pages).

\re{[14]} T. Kato, Perturbation Theory for Linear Operators, 2nd ed.,
Springer-Verlag, Berlin \newline/Heidelberg/New York/Tokyo, 1984.

\re{[15]} S. J. Lee, M. Z. Nashed, Least-squares solutions of multi-valued linear
operator equations in Hilbert spaces, J. Approx. Theory 38(1983) 380--391.

\re{[16]} S. J. Lee, M. Z. Nashed, Constrained least-squares solutions of linear
inclusions and singular control problems in Hilbert spaces, Appl. Math.
Optim. 19(1989) 225--242.

\re{[17]} S. J. Lee, M. Z. Nashed, Algebraic and topological selections of
multi-valued linear relations, Ann. Scuola Norm. Sup. Piss 17(1990) 111--126.

\re{[18]} S. J. Lee, M. Z. Nashed, Normed linear relations:
domain decomposability, adjoint subspaces and selections,
Linear Algebra Appl. 153(1991) 135--159.

\re{[19]} M. Lesch, M. Malamud,
On the deficiency indices and self-adjointness of
symmetric Hamiltonian systems, J. Differ. Equations 18(2003) 556--615.

\re{[20]} Y. Liu, Y. Shi, Regular approximations of spectra of singular second-order
   symmetric linear difference equations, Linear Algebra Appl. 444(2014) 183--210.

\re{[21]} D. G. Luenberger, Optimization by Vector Space Methods, Wiley, New
York, 1969.

\re{[22]} M. Z. Nashed, Operator parts and generalized inverses of linear manifolds
with applications, in Trends in Theory and Practice of Nonlinear Differential Equutions
(V. Lakshmikantham, Ed.), Marcel Dekker, New York, 1984, pp. 395--412.

\re{[23]} J. von Neumann, \"Uber adjungierte Funktional-operatoren,
Ann. Math. 33(1932) 294--340.

\re{[24]} J. von Neumann, Functtional Operators, Vol. II: the Geometry of Orthogonal
Spaces, Ann. of Math. Stud. 22, Princeton U.P., 1950.

\re{[25]} M. Reed, B. Simon, Methods of Modern Mathematical Physics
II: Fourier Analysis, Self-adjointness, Academic Press, 1972.

\re{[26]} G. Ren, Y. Shi,
Defect indices and definiteness conditions
for discrete linear Hamiltonian systems, Appl. Math. Comput.
218(2011) 3414--3429.

\re{[27]}
G. Ren, Y. Shi, Self-adjoint extensions of for a class of discrete linear
Hamiltonian systems, Linear Algebra Appl. 454 (2014) 1--48.

\re{[28]} S. M. Robinson, Normed convex processes, Trans. Amer. Math. Soc.
174(1972) 124--140.

\re{[29]} K. Schm\"udgen, Unbounded Self-adjoint Operators on Hilbert space,
Springer Dordrecht Heidelberg/ New York/ London, 2012.

\re{[30]} Y. Shi, The Glazman-Krein-Naimark theory for Hermitian
subspaces, J. Operat. Theor. 68(1)(2012) 241--256.

\re{[31]} Y. Shi, Stability of essential spectra of self-adjoint subspaces under compact
   perturbations. J. Math. Anal. Appl. 433(2016) 832-851.

\re{[32]} Y. Shi, C. Shao, Y. Liu,
Resolvent convergence and spectral approximations
of sequences of self-adjoint subspaces, J. Math. Anal. Appl. 409(2014) 1005--1020.

\re{[33]} Y. Shi, C. Shao, G. Ren, Spectral properties of self-adjoint subspaces,
Linear Algebra Appl. 438(2013) 191--218.

\re{[34]} Y. Shi, H. Sun,
Self-adjoint extensions for second-order
symmetric linear difference equations, Linear Algebra Appl.
434(2011) 903--930.

\re{[35]} H. Sun, Y. Shi, Spectral properties of singular discrete linear Hamiltonian systems,
J. Diff. Equ. Appl. 20(2014) 379--405.

\re{[36]} J. Weidmann, Linear Operators in Hilbert Spaces,
Graduate Texts in Math., vol.68, Springer-Verlag, New
York/Berlin/Heidelberg/Tokyo, 1980.

\re{[37]} D. Wilcox, Essential spectra of linear relations,
Linear Algebra Appl. 462(2014) 110--125.

\end{document}